\documentclass[11pt]{amsart}   
\usepackage{amsmath}
\usepackage{amsfonts}
\usepackage{amssymb}
\usepackage{latexsym}
\usepackage{graphicx}%

\newtheorem{theorem}{Theorem}[section]
\newtheorem{lemma}[theorem]{Lemma}

\newtheorem{corollary}[theorem]{Corollary}
\theoremstyle{definition}

\theoremstyle{remark}

\numberwithin{equation}{section}

\def\Re{\mathbb{R}}

\def\Fi{\mathcal{F}}

\begin{document}

\vspace{0.5in}

\title[]%
{Continuity of the It\^o-map for H\"older rough paths with applications to the Support Theorem in H\"older norm}

\date{\today}

\author{Peter K. Friz}
\address{Courant Institute, New York University}
\email{Peter.Friz@cims.nyu.edu}




\begin{abstract} 
Rough Path theory is currently formulated in $p$-variation topology.
We show that in the context of Brownian motion, enhanced to a Rough
 Path, a more natural H\"older metric $\rho$ can be used.
 Based on fine-estimates in Lyons' celebrated Universal Limit Theorem we obtain Lipschitz-continuity of the It\^o-map (between Rough Path spaces equipped with $\rho$).
We then consider a number of approximations to Brownian Rough Paths and establish their convergence w.r.t. $\rho$. 
In combination with our H\"older ULT this allows sharper results than the
$p$-variation theory. Also, our formulation avoids
the so-called control functions and may be easier to use for
non Rough Path specialists.
As concrete application, we combine our results with 
ideas from ~\cite{MS} and ~\cite{LQZ} and obtain the Stroock-Varadhan Support Theorem in
H\"older topology as immediate corollary.
\end{abstract}

\maketitle

\section{Introduction}

\subsection{Background in Rough Path theory} 

Over the last years T.Lyons and co-authors developed a general theory
of integration and differential equations of form
\begin{equation}
      dy_t = f(y_t) dx_t  \label{ODE}.
\end{equation}
To include the important example of stochastic differential equations,
$x$ must to be allowed to be "rough" in some sense.  Standard H\"older-regularity of Brownian motion, for instance,
 implies finite $p$-variation only for $p>2$.
Another issue is to explain (deterministically) the difference between SDEs based on Stratonovich- vs. It\^o-integrals.
Last not least, motivated from examples like Fractional Brownian motion, 
driving signals much rougher than Brownian motion should be included. \\
All this has been accomplished in a very satisfying way and the reader can nowadays find the general theory
exposed in ~\cite{L98}, ~\cite{LQ}, ~\cite{Le}.

 Loosely speaking, for general $p \ge 1$, one needs to "enhance" the driving signal $x$, with values in some Banach-space $V$,
to $X \in  V \oplus V^{\otimes 2} ... \oplus V^{\otimes [p]}$
such that the resulting object $X$
satisfies certain algebraic \footnote{For {\it algebraic} convenience $X$
is often enhanced to $\Re \oplus V \oplus V^{\otimes 2} ... \oplus V^{\otimes [p]}$ with scalar component constant $1$.} and analytic
conditions. For $x$ of finite variation, this enhancement will simply consist of all the
iterated integrals of $x$,
$$
      X^k_{s,t} := \int_{s < u_1 < ... < u_k < t} dx_{u_1} \otimes ...
						\otimes dx_{u_k}, \ \ \
						k = 1, ... , [p].
$$
 These are the {\it Smooth Rough Paths.} Consider a time-horizon of $[0,1]$
 (valid for the rest of the paper) and introduce the
$p$-variation metric, defined as
$$
     d(X,Y) = \max_{k=1,...,[p]} 
             \Big( \sup_D \sum_l | X^k_{t_{l-1},t_l} - Y^k_{t_{l-1},t_l} |^{p/k} 		\Big)^{k/p},
$$
where $\sup_D$ runs over all finite divisions of $[0,1]$. Here $|.|$ denotes (compatible) tensor norms in $V^{\otimes k}$. Closure of Smooth Rough Paths w.r.t.
this metric yields the class of {\it Geometric Rough Paths}
\footnote{A slightly weaker definition of Geometric Rough Path appeared
in ~\cite{L98}. We follow the more recent ~\cite{LQ}.}, denoted by
 $G\Omega_p(V)$.
The solution-map, also called {\it It\^o-map}, to (\ref{ODE}) is then a continuous map from $G\Omega_p(V) \to G\Omega_p (W)$, provided
$f: W \to L(V,W)$ satisfies mild regularity conditions. This is Lyons'
celebrated {\it Universal Limit Theorem}. In particular,
smooth approximations $X(n)$ which converge in $p$-variation to $X \in G\Omega_p(V)$ will cause the corresponding solutions $Y(n)$ to converge to $Y$ in
$p$-variation. Hence, one deals with some kind of generalized Stratonovich theory. 

{\it However}, the so important case of $p\in (2,3)$, on which this
paper will focus, allows for more. Following ~\cite{LQ} p149 and
 also ~\cite{L98} the driving signal only needs to be a {\it Rough Path of finite $p$-variation}. By definition, this is a continuous map
$$
   (s,t) \to (X^1_{s,t},X^2_{s,t}) \in V \oplus V^{\otimes 2},
$$
where $0 \le s \le t \le 1$, satisfying the algebraic condition
\begin{equation} 
 X^1_{s,u} = X^1_{s,t} + X^1_{t,u}
\ \ \mbox{and} \ \ \
 X^2_{s,u} = X^2_{s,t} + X^2_{t,u} + X^1_{s,t} \otimes X^1_{t,u},
\label{chen}
\end{equation}
whenever $s \le t \le u$, and the analytic condition $d(X,0) < \infty$ i.e.
\begin{equation}
     \sup_D \sum_l | X^k_{t_{l-1},t_l} |^{p/k} < \infty, \ \ \ \ k =1,2 \ \ \ \label{finite_p_var}.
\end{equation}
(Sometimes we refer to $k=1,2$ as {\it first} resp. {\it second level}.)
The class of Rough Paths of finite $p$-variation is denoted by $\Omega_p(V)$.
 Clearly,
$$
    \{ \mbox{Smooth Rough Paths} \} \subset G\Omega_p(V) \subset \Omega_p (V).
$$ 
Condition (\ref{chen}) is known as {\it Chen relation} and expresses simple
additive properties whenever $X^2$ is obtained as {\it some} iterated integral
including the cases of Stratonovich- resp. It\^o-Enhanced Brownian Motion (EBM), where 
 $$
X^1_{s,t} = \beta^1_{s,t} = \omega_t - \omega_s,
$$
the increments of Brownian Motion on $V=\Re^d$, defined on the
usual Wiener-space $(C([0,1],\Re^d), \mu)$.
$X^2_{s,t} =
\beta^2_{s,t}$ is defined via stochastic integration,
\begin{equation}
    \int_s^t (\omega^i_r - \omega^i_s ) \circ d\omega^j_r 
   \ \ \mbox{resp.} \ \  \int_s^t (\omega^i_r - \omega^i_s ) d\omega^j_r.
     \label{2nd_level}
\end{equation}
Clearly, there is a modification of $\beta^2$, denoted by the same name,
such that $\beta^2$ (and hence $\beta$) is continuous in $s,t$ for all $\omega$. With this choice, it is well-known that (both Stratonovich- and It\^o-)EBM are a.s. elements 
of $\Omega_p(V)$, also called {\it Brownian Rough Paths}.  Indeed, this will follow from a.e. Stratonovich-EBM $\beta = \beta(\omega) = (\beta^1,\beta^2)$ being a Geometric Rough Path and this is the (somewhat computational) approach
 in ~\cite{LQ}, ~\cite{Le}.

We emphasize once more, that for $p \in (2,3)$ and driving signal
in $\Omega_p (V)$ a complete theory is available, covering
both It\^o- and Stratonovich-SDEs. In particular, the Universal Limit
 Theorem holds, see ~\cite{LQ}, p164 and Section \ref{sec4}. \\

\subsection{Definitions and outline}

Our focus being on Brownian Rough Paths, let $p\in (2,3)$ and $V=\Re^d$
from here on.
Given a Rough Path
\begin{equation}
   \| X \| := \max_{k=1,2} \sup_{0 \le s < t \le 1} \frac{ |X_{s,t}^k| }{ |t-s|^{k/p} }.     \label{Hnorm}
\end{equation}
We call Rough Paths with $\|X \| < \infty$ {\it H\"older Rough Paths} and
write $X \in H\Omega_p(V)$. This condition is similar to (\ref{finite_p_var})
but stronger. Note $H\Omega_p(V) \subset \Omega_p(V)$. Clearly, $H\Omega_p(V)$
becomes a metric space under $\rho(X,Y) := \| X-Y \|$.
(The norm $\| . \|$  appeared in ~\cite{L95} but seems to have disappeared
in the current {\it $p$-variation } Rough Path theory.) 
In {\bf Section \ref{sec2}} we use Kolmogorov's criterion to prove
a 2-parameter H\"older-regularity of an object intimately related
to $\beta$.  This implies $\beta \in H\Omega_p(V)$.
(Such a result was mentioned without proof in ~\cite{L95}
and might be contained in unpublished thesis-work, ~\cite{S}.) As further corollary, we have the well-known $\beta \in \Omega_p(V)$. 
In {\bf Section \ref{sec3}}, we establish convergence in $\rho$-metric
of piecewise linear approximations based on nested partitions. This gives
a simple and novel proof that Stratonovich-EBM is an element of $G\Omega_p (V)$.(Our proof combines the preceding regularity result with
soft martingale and compactness arguments. In a sense, all dyadic
 approximations done by hand in ~\cite{LQ}, ~\cite{Le} are isolated in our initial application of Kolmogorov's criterion.) At last, the results are extended
to adapted dyadics approximations.

{\bf Section \ref{sec4}} is a recall of Lyons' Universal Limit Theorem.
In {\bf Section \ref{sec5}} 
we translate his fine estimates in terms of {\it control functions} to
plein Lipschitz-continuity of the It\^o-map
from $(H\Omega_p(V),\rho)$ to $(H\Omega_p(W),\rho)$.

Finally, in {\bf Section \ref{sec6}} we apply our refined Rough Path machinery to prove the celebrated Support Theorem in H\"older-topology.
 As observed in ~\cite{MS}, the proof can be reduced to two convergence results (one for each inclusion) and these follow immediately from our results. 

{\bf Remarks:}
-The use of Rough Path theory to prove the Support Theorem was first carried out in ~\cite{LQZ}.  Using the standard (= $p$-variation) Rough Path machinery the well-known H\"older-topology result ( ~\cite{BGL}, ~\cite{MS}, ~\cite{ST}) was not recovered. We also note, that the approach in ~\cite{LQZ} relies on correlation inequalities.

-A recent preprint, ~\cite{Gu}, proposes to re-prove rough path theory in a
pure H\"older context, apparently motivated by (\ref{Hnorm}).

-Forthcoming joint-work with Nicolas Victoir will (among other things) contain a rough paths proof of the support theorem in yet a stronger norm (modulus norm).\footnote{The overlap to ~\cite{FV} is kept to a minimum. In particular, 2-paramter regularity and convergence, definition (\ref{Hnorm}) and Lipschitz-continuity of It\^o-map are not contained in ~\cite{FV}.}

\section{H\"older-regularity of Enhanced Brownian motion} \label{sec2}

As before, let $(\beta^1, \beta^2) = \beta = \beta(\omega)$ be (It\^o -or Stratonovich) EBM based on $d$-dimensional BM, chosen s.t. $(s,t) \mapsto
\beta_{s,t} (\omega)$ is continuous for all $\omega$.
Recall $0 \le s,t \le 1, \ p \in(2,3)$.

\begin{theorem} \label{thm21}
Let $\gamma \in (0,1/2 - 1/p)$.
 Introduce the following $2$-parameter processes 
\begin{equation}
                    Z^k_{s,t} = \frac{\beta^k_{s,t}}{|t-s|^{k/p}}
                     \in V^{\otimes k} \cong \Re^{d^k}  \ \ \   k=1,2,
                     \label{defZ}
\end{equation}
whenever $t>s$ and set them zero otherwise. Then for $\mu$-a.e. $\omega$
$Z^k=Z^k(\omega)$, as function of $s,t$, is H\"older continuous of any exponent $k\gamma$. Also,
\begin{equation}
          \sup_{s \ne s', t \ne t'} \frac{|Z^k_{s,t} - Z^k_{s',t'}|}
                             { [|t'-t| + |s'-s|]^{k\gamma} }
          \ \ \in L^{\infty-} := \cap_{q<\infty} L^q.
          \label{HregZ}
\end{equation}
\end{theorem}

The following lemma is proved in the appendix.
\begin{lemma} \label{lemma21}
Fix $i,j \in \{ 1,...,d \}$ and set $Z=Z^{2,ij}$. Then there is $c=c(p)$ s.t.
$$
 E[ | Z_{s,t} - Z_{s',t'} |^2] \le c [|s'-s| + |t'-t|]^{2(1-2/p)}.
$$
\end{lemma}

{\bf Proof of Theorem \ref{thm21}.} 

We only prove $k=2$. The first level i.e. $k=1$ is similar but easier.
As before, fix $i,j$ and set $Z_{s,t} = Z^{k,ij}_{s,t}$. By equivalence of
finite dimensional norms, it suffices to show that $Z$ is H\"older.
We first consider It\^o-EBM.
With this choice $Z_{s,t}$ is an element of the second Wiener-It\^o-chaos
and a moment estimate similar to the Gaussian case is available (~\cite{RY}, p207). With lemma \ref{lemma21},
\begin{eqnarray*}
 E[ | Z_{s,t} - Z_{s',t'} |^q ] & \le & c(q) \Big(  E[ | Z_{s,t} - Z_{s',t'} |^2 ] \Big)^{q/2} \\
                                & \le & c [|s'-s| + |t'-t|]^{(1- 2/p)q}.
\end{eqnarray*}
We can choose $q$ arbitrarily large and, by Kolmogorov's criterion as found
in ~\cite{RY}, obtain H\"older regularity
for any exponent $\gamma$ less than $1-2/p$. \\
The only thing left to consider is $Z=Z^{2,ij}$ based on Stratonovich-EBM. 
Only on the diagonal $i=j$ there is a non-zero difference (coming
from the quadratic variation of BM),
$$
      Z_{s,t}^{\mbox{\tiny It\^o}} - Z_{s,t}^{\mbox{\tiny Stratonovich}}
      = \frac{1}{2}|t-s|^{1-2/p}.
$$
As composition of the Lipschitz-map $(s,t) \to |t-s|$ and a $(1-2/p)$-H\"older
map, the map $(s,t) \mapsto |t-s|^{1-2/p}$ is itself $(1-2/p)$-H\"older.
This implies the $\gamma$-H\"older-regularity of $Z=Z^{\mbox{\tiny Stratonovich}}$. The proof is finished. QED

\begin{corollary} \label{cor1}
Enhanced Brownian motion $\beta$ is $\mu$-a.s. an element of $H\Omega_p(V) \subset \Omega_p(V)$.  More precisely, there exists $C=C(\omega) \in L^{\infty-}$ such that
for all $s \le t$ and $k=1,2$ 
\begin{equation}
          | \beta^k_{s,t} |^{p/k} \le C |t-s|.    \label{control}
\end{equation}
(In Lyons' terminology, $\beta$ admits the
additive control $C|t-s|$.)
\end{corollary}

{\bf Proof of Corollary \ref{cor1}.} 

The continuous function $(s,t) \mapsto Z^k_{s,t}$ achieves its maximum
$C(\omega)$ which is easily estimated by the (\ref{HregZ}).   QED

{\bf Remark: } We showed that a.s. $Z^k_{..}$ is an element of the
H\"older-space $C^{0,k\gamma}([0,1]^2,V^{\otimes k})$. Moreover,
(\ref{HregZ}) can be used as a norm.

\section{Approximations to Brownian Rough Paths}   \label{sec3}

\subsection{Piecewise linear nested approximations}

Based on piecewise linear nested approximations $\omega(n)$ of the underlying Brownian path $\omega$ we construct a Smooth Rough Path denoted by 
$\beta(n) = (\beta^1(n),\beta^2(n))$.
We assume that these partitions are deterministically chosen and
that their mesh goes to $0$. Note that the commonly used piecewise linear
{\it dyadic} approximations
(see ~\cite{ST}, ~\cite{M}, ~\cite{LQ} ,...) fall into the considered class.

As before, introduce $Z^1, Z^2$ based on Stratonovich-EBM $\beta$
and, similarly $Z^1(n), Z^2(n)$ based on $\beta(n)$.
 The matrix-valued processes $Z^2, Z^2(n), ...$ split into
symmetric and anti-symmetric parts. For instance,
$$
       Z^2 = \hat Z^2 + A^2,
$$
where $\hat .$ indicates symmetrization.

\begin{theorem} \label{thm31}
 Let $p \in (2,3), \gamma \in (0,1/2-1/p)$. For first level and symmetric part
of the second level, convergence of the approximations holds a.s.
in H\"older-space of exponent $\gamma$,
\begin{eqnarray}
              \| Z^1 - Z^1 (n) \|_{C^{0,\gamma}} \to 0,  \label{l1} \\
              \| \hat Z^2 - \hat Z^2 (n) \|_{C^{0,\gamma}} \to 0. 
							\label{l2s} 
\end{eqnarray}
 For the antisymmeric part $A^2$, based on the L\'evy-area of the
the underlying BM, we have a.s.
\begin{equation}
              \| A^2 - A^2 (n) \|_{C^{0,2\gamma}} \to 0.    \label{l2a}
\end{equation}
\end{theorem}

{\bf Remark: }Given that $Z^2$ is itself $2\gamma$-H\"older it is quite possible
that $Z^2(n) \to Z^2$ in $C^{0,2\gamma}$, but none of our conclusions will rely on this.

\begin{corollary} \label{cor32}
 There exists a random constant $C < \infty$ a.s.
and a random sequence $a_n(\omega) \to 0$ a.s.
such that, for $s \le t$, $k=1,2$ 
\begin{equation}
      | \beta_{s,t}^k (n) |, | \beta_{s,t}^k |
			 \le [C(\omega) |t-s|]^{k/p}   \label{unifest}
\end{equation}
and
\begin{equation}
      | \beta_{s,t}^k - \beta_{s,t}^k (n) | \le a_n (\omega) |t-s|^{k/p}.
      \label{randomsequ}
\end{equation}
(In Lyons' terminology, the additive control $C(\omega) |t-s|$ is uniform
for the entire sequence $\beta(n)$ and controls the convergence.  \footnote{Strictly speaking, $(C \lor 1) |t-s|$ will be the required control.} )
\end{corollary}

These estimates translate to

\begin{corollary} \label{cor32prime}
The (Smooth Rough Path-) approximations $\beta(n)$ converge to $\beta$ in H\"older-metric $\rho$,
$$
	\rho ( \beta(n), \beta ) \to 0  \ \ \mbox{a.s.}
$$
(Since this implies convergence in $p$-variation metric $d$ we identify,
en passant, Stratonovich-EBM $\beta$ as a Geometric Rough Path.)
\end{corollary}


{\bf Proof of Theorem \ref{thm31}.}  We are able to do most
of the work for levels $k=1,2$ at the same time. For the moment,
fix $i \ne j$. For $k=1$ set 
$$
      Z = Z^{1,i},
$$
while for $k=2$ set
$$ 
     Z= Z^{2,ij}.
$$
Either way, we have a real-valued $2$-paramter process, 
$k\gamma$-H\"older according to Theorem \ref{thm21}. This means
for some $L^{\infty-}$-r.v. $L$ we have the inequality
$$
-L [ |t' -t| + |s' - s| ]^{k\gamma} \le
Z_{s,t} - Z_{s',t'} \le
L [ |t' -t| + |s' - s| ]^{k\gamma}.
$$
Now condition w.r.t. $\Fi_n := \sigma( \beta_{k/2^n}:
 k = 0,..., 2^n )$. Set $L_n = E[ L | \Fi_n]$. 
It is not hard to see (~\cite{M}, p216)
 that  $E [ Z | \Fi_n ] = Z(n)$.
(In the case $k=2$, here is where we use $i \ne j$.) Hence,
$$
-L_n [ |t' -t| + |s' - s| ]^{k\gamma} \le
Z_{s,t} (n) - Z_{s',t'} (n) \le
L_n [ |t' -t| + |s' - s| ]^{k\gamma}.
$$
Note that $L_n$ is an $L^{\infty-}$-bounded martingale.
By Doob's $L^p$-inequality we see that,
$$
     C(\omega) := \sup_n L_n (\omega)
$$
is also in $L^{\infty-}$. Consequently,
$$
    \sup_{s,t,s',t'} \frac{|Z_{s,t}(n) - Z_{s',t'}(n)|}{[ |t' -t| + |s' - s| ]^{k\gamma}} \le L_n \le C.
$$
This implies that the sequence $(Z(n))$ is bounded in the
 H\"older-space $C^{0,k\gamma}$. We could have started with
$\gamma + \epsilon$ as long as $\gamma + \epsilon < 1/2 - 1/p$.
Then the conclusion is boundedness in $C^{0,k(\gamma+\epsilon)}$
and by compactness there is a convergent subsequence in
$C^{0,k\gamma}$. But every possible limit point is identified
as $Z$ itself, since for $s,t$ fixed, $Z_{s,t}(n) \to Z_{s,t}$
by martingale convergence. This implies that $Z(n)$ actually
converges to $Z$ in $C^{0,k\gamma}$ and the same holds true
for the antisymmetric part of $Z^2_{s,t}$ itself since
all the diagonal, $i=j$, is zero. 

At last, we need to consider the case
$$
   Y_{s,t} := 2Z^{2,ii}_{s,t} = \frac{ (\beta^i_t - \beta^i_s)^2}{ |t-s|^{2/p} } = (Z^{1,i}_{s,t})^2 =: (Z_{s,t})^2.
$$
Similarly define  $Y(n), Z(n)$.
We claim that $Y(n)$ tends to $Y$ in $C^{0,\gamma}$. Certainly,
for $s,t$ fixed $Y_{s,t} (n) \to Y_{s,t}$ which identifies
every possible limit point of $Y(n)$. Hence it suffices,
by the same compactness argument as before, to show that $Y(n)$ is uniformly bounded in $C^{0,\gamma}$. But $Y(n) = [Z(n)]^2$ with
$Z(n) = Z^{1,i} (n)$ and this last sequence was shown earlier to be 
uniformly bounded (even convergent) in $C^{0,\gamma}$.
On the other hand, the map $x \mapsto x^2$ is (locally) Lipschitz,
and since $\{Z(n)\}$ remains in a ball in $C^{0,\gamma}$ we conclude
that $\{Y(n) \}$ remains in a (possibly larger) ball in $C^{0,\gamma}$ 
as well. The claim is proved. \\
Together with the earlier results for $k=2, i \ne j$ we find
that $\hat Z^2 (n)$, the symmetric part of $Z^2(n)$, converges in
$C^{0,\gamma}$ to $\hat Z$. The proof is finished.   QED

{\bf Proof of Corollary \ref{cor32prime}.}

For $k=1$ the estimates
 (\ref{unifest}) and (\ref{randomsequ}) are
 an immediate consequence of
(\ref{l1}) and $Z^1 \in C^{0,\gamma}$. 


  Towards $k=2$ we first
consider the anti-symmetric part
of $\beta^2_{s,t}  - \beta^2_{s,t} (n)$. But this is simply bounded by
$$
     a_1 (n,\omega) |t-s|^{2/p}
$$
where $a_1 (n)$ denotes the l.h.s. of (\ref{l2a}). Towards the
symmetric part of $\beta^2_{s,t}  - \beta^2_{s,t} (n)$ observe that
$$
      \hat \beta^{2,ij}_{s,t} = \frac{1}{2}\beta^{1,i}_{s,t} \beta^{1,j}_{s,t}.
$$
It suffices to estimate one component of $\hat \beta^2  - \hat \beta^2(n)$,
namely,
\begin{eqnarray*}
    | \beta^{1,i}_{s,t} \beta^{1,j}_{s,t}
     -     \beta^{1,i}_{s,t} (n) \beta^{1,j}_{s,t} (n) |
    & \le & |\beta^{1,i}_{s,t}| |\beta^{1,j}_{s,t} -\beta^{1,j}_{s,t}(n) | \\
    & & +|\beta^{1,i}_{s,t} -\beta^{1,i}_{s,t}(n) | | \beta^{1,j}_{s,t}(n) | \\
\end{eqnarray*}
for arbitrary $i,j$. From (\ref{l1}) it follows that
$$
     |\beta^1_{s,t} - \beta^1_{s,t} (n) | \le a_2(n,\omega) |t-s|^{1/p},
$$
where $a_2 (n)$ denotes the l.h.s. of (\ref{l1}). Together with
the uniform estimates (\ref{unifest}) we conclude that
$$
    |\hat \beta^2_{s,t}  - \hat \beta^2_{s,t}(n)| \le a_3 (n,\omega) |t-s|^{2/p},
$$
where $a_3$ is a deterministic constant times $C(\omega)a_2$. 
For
$$
     a(n) := \max \{ a_1 (n), a_3 (n) \}
$$
estimate (\ref{randomsequ}) will then hold true. 
Also, (\ref{unifest}) follows by the triangle-inequality 
and the regularity of $\beta^2$. QED.

\subsection{Adapted dyadic approximations}

Now let $\omega(n)$ be the dyadic piecewise linear approximation to
 a Brownian path $\omega$, i.e. piecewise linear from $\omega_{i/2^n}$ to
 $\omega_{(i+1)/2^n}$.  Note, that $\omega(n)$ is not adapted to the
 Brownian filtration. This
suggests to look at the following {\it adapted} approximation,
$$
      \omega^{ad}_t (n) = \omega_{(t - 2^{-n}) \lor 0} (n).
$$
We can lift the path $\omega^{ad}(n)$ to make it a (Smooth) Rough Path,
which we denote by $\beta^{ad}(n)$. Similarly, and as before, $\omega(n)$ is lifted to $\beta(n)$.

\begin{corollary} \label{cor34}
The Smooth Rough Paths  $\beta^{ad} (n)$ converge a.s. to
Stratonovich-EBM $\beta$ in H\"older-metric
$\rho$ (and consequently in $p$-variation metric $d$).
\end{corollary}

{\bf Proof: } Introduce a shift-operator on the path-level
s.t. for any path $x(t) \in V$,

$$
	\tau^\epsilon: x(.)  \mapsto x((. - \epsilon) \lor 0).
$$
This lifts to a map on Rough Paths. With $\epsilon = 1/2^n$,
Then,
$$
       \beta^{ad} (n) = \tau^\epsilon ( \beta (n) ).
$$
In particular, for the second level,
$$
       \beta^{ad,2}_{s,t} (n) = \beta^{2}_{(s-\epsilon) \lor 0, 
       		(t-\epsilon) \lor 0}.
$$

By the triangle inequality,
$$
   \rho ( \tau^\epsilon( \beta (n) ), \beta )
   \le
   \rho (\tau^\epsilon( \beta (n) ) , \tau^\epsilon( \beta ) )
   +
   \rho (\tau^\epsilon( \beta ) , \beta ).
$$
Since $\rho(\beta(n),\beta) \to 0$ so does the first term on the r.h.s.
($\rho$ is insensitive to shift). So all that remains to show
is that $ \rho (\tau^\epsilon( \beta ) , \beta ) \to 0$ as $\epsilon$ tends
to zero. Written out, this means 
$$
\max_{k=1,2} \sup_{0 \le s < t \le 1} \frac{ |\beta_{(s-\epsilon)\lor 0,(t-\epsilon) \lor 0}^k 
- \beta_{s,t}^k| 
}{ |t-s|^{k/p} }
$$
goes to zero with $\epsilon$. We can estimate this by
$$
   \sup_{0 \le s < t \le \epsilon} ( ...) +
   \sup_{0 \le s < \epsilon \le t \le 1} (...) +
   \sup_{\epsilon \le s < t \le 1} + (...).
$$
Each part is easily seen to converge to $0$ with $\epsilon$ by using
the $(s,t)$ H\"older-property of $Z^k_{s,t} = \beta^k_{s,t} / |t-s|^{k/p}$
established in Theorem \ref{thm31}.   QED

Clearly, the last corollary implies that on path level and
in H\"older-norm with exponent less than $1/p$,
$$
      \omega^{ad} (n) \to \omega.
$$
Then, trivially,
$$
      T_n (\omega) = \omega - \omega^{ad}(n) + h \to h,
$$
for, say, any piecewise linear dyadic path $h$.
 We will lift this convergence result to Rough Path level.
As in ~\cite{LQZ} we shall denote the lifts of dyadic piecewise linear paths by  $\mathcal D$, a set
of smooth rough paths. Before doing so,
recall the Stratonovich enhancement of Brownian motion,
$$
     \omega \to \beta_{s,t} (\omega) = ( \omega_t - \omega_s,
                \frac{1}{2} (\omega_t - \omega_s)^{\otimes 2} + A_{s,t} )
                 \in V \oplus V^{\otimes 2}.
$$
Due to the L\'evy-area $A_{s,t}=A_{s,t}(\omega)$ this is only an a.s. defined
function of $\omega$  (although we picked a modification, determined
up to indistinguishability, s.t. $(s,t) \to
\beta_{s,t}(\omega)$ is continuous). We saw in section \ref{sec2} that 
$\beta \in H\Omega_p(V)$ for $\mu$-a.e. $\omega$.
By Girsanov's theorem, $\beta (T_n (\omega)$ is well-defined and
in $H\Omega_p(V)$ for $\mu$-a.e. $\omega$. 

\begin{corollary} \label{cor36}
The H\"older Rough Paths  $ \beta( T^h_n (\omega)) $ converge a.s. in 
H\"older-metric $\rho$ to the (Smooth) Rough Path $(h^1,h^2) \in {\mathcal D}$, associated to the piecewise linear dyadic path $h(.)$.
\end{corollary}

{\bf Proof:} It suffices to consider the L\'evy-area, more specifically off-diagonal term
of the 2nd-level.
That is, we want to show that, for $i \ne j$,
$$
   \sup_{s<t} \frac{|\beta^{2,ij}_{s,t} (T^h_n (\omega)) - h^{2,ij}_{s,t} |}
                   { |t-s|^{2/p} }
$$
tends to zero.
To this end, a Riemann-sum approximation shows that the following expansion holds a.s.
(we omit $s,t$, the following integrations are understood over the simplex
$\{ (u_1,u_2): s \le u_1 \le u_2 \le t \}$).
\begin{eqnarray*} 
&  &	\beta^{2,ij} (T^h_n (\omega))  =  \beta^{2,ij} (\omega - \omega^{ad}(n)+h)  =  h^{2,ij} \\
   & + & \beta^{2,ij} (\omega) 
    +     \int d \omega^{ad,i}(n) d\omega^{ad,j}(n) 
     - \int d \omega^i d\omega^{ad,j}(n)
	-  \int d \omega^{ad,i}(n) d\omega^j  \\
   & + & \int d h^i d\omega^j - \int d h^i d\omega^{ad,j}(n)  
    +  \int d \omega^i dh^j - \int d \omega^{ad,i}(n) dh^j.
\end{eqnarray*}
(All iterated integrals here make sense as Young-integrals.)
 Note that the last corollary implies that
\begin{equation}
    \int d \omega^{ad,i}(n) d\omega^{ad,j}(n)
    \to 
    \beta^{2,ij} (\omega)    \label{ad_conv}
\end{equation}
w.r.t. the second level part of the metric $\rho$ (just keep $k=2$ in
its definition).
Clearly, $\int d \omega^i d\omega^{ad,j}(n)$ and
	$\int d \omega^{ad,i}(n) d\omega^j$ are even better
approximations to $\beta^{2,ij}(\omega)$ and hence converge to the latter.  \footnote{This
is easily made precise by partial conditioning in the proof of Theorem \ref{thm31}, that is condition w.r.t. the $i$th resp. $j$th component only.}
Then all four terms together in the second line (r.h.s.) above will
converge to zero. \\
 As for the third line, observe that (\ref{ad_conv}) holds
for a.e. Brownian path $\omega$. By Girsanov's theorem we have a.s. convergence
after replacing $\omega$ by $\omega + h$. Now use an expansion as before,
$$
    \beta^{2,ij}(\omega+h) = \beta^{2,ij}(\omega) + h^{2,ij}
    +  \int dh^i d\omega^j +  \int d\omega^i dh^j.
$$
Similarly, expand 
 $$
\int d(\omega^{ad,i}(n) + h^i(n)) d(\omega^{ad,j}(n) + h^j(n)).
$$
(Note that $h(n) = h \in {\mathcal D}$ for all $n$ large enough.)
 As already mentioned,
(\ref{ad_conv}) still holds after replacing $\omega$ by $\omega+h$.
Together with the expansions, this gives exactly the required cancelation
(as $n\to \infty)$ of the third line above.    QED

{\bf Remark:} The preceding proof involves a perturbation of the rough
path $\beta$ (essentially) in a Cameron-Martin-type direction, $(h^1,h^2)$.  Such and more
 general perturbations have been studied systematically in ~\cite{LQ97}, see also ~\cite{LQ}.
 Indeed, we could have a based our proof on some of their general results.


\section{A primer on the Universal Limit Theorem}    \label{sec4}

We just summarize and plug together a few statements from ~\cite{LQ}. 
Set $V=\Re^d, W=\Re^N$. As before, $p \in (2,3)$.

 Recall that
a {\it control function}\footnote{The reader who does not want to
know about control functions may jump to the statement of Theorem
\ref{thm52} directly.} is, by definition, a
 non-negative continuous function $w$ on $\{ 0 \le s \le t \le 1 \}$,
super-additive, that is, $w(s,t) + w(t,u) \le w(s,u)$ and hence zero on the diagonal.

\begin{theorem} \label{ULT}
Let $f \in C^3 (W,L(V,W))$, globally Lipschitz. Then the It\^o-map 
from $\Omega_p (V) \to \Omega_p (W)$, obtained from ``enhancing''
$$
     dy = f(y)dx,	y(0) = y_0,
$$
 is continuous w.r.t. the $p$-variation metric.
Moreover, if $w=w(s,t)$ is a control function
such that (always for all $s,t$)
\begin{equation}
     | X^k_{s,t} |, | \hat X^k_{s,t} | \le w(s,t)^{k/p}, \ \ \ k=1,2    \label{ass1}
\end{equation}
and
\begin{equation}
   | X^k_{s,t} - \hat X^k_{s,t} | \le \epsilon w(s,t)^{k/p}, \ \ \ k=1,2,    \label{ass2}
\end{equation}
then there exists a constant $C=C( \max w)=C(\max w, f,p,y_0)$ such that
$$
     | Y^k_{s,t} |, | \hat Y^k_{s,t} | \le (C w(s,t))^{k/p}, \ \ \ k=1,2
$$
and
$$
   | Y^k_{s,t} - \hat Y^k_{s,t} | \le \epsilon (C w(s,t))^{k/p}, \ \ \ k=1,2.
$$
\end{theorem}

{\bf Proof: } The fine estimate in terms of control functions $w$ are stated and proved in ~\cite{LQ}, p163.
For the readers convenience, let us quickly show how to obtain continuity w.r.t $p$-variation metric $d(.,.)$.
Assume $d(X(n),X) \to 0$. Clearly, for any subquence, $d(X(n_i), X) \to 0$. By
\footnote{Unfortunately, the cited page contains a misprint.
The supremum of superadditive functions may loose this property. However,
all this is readily fixed and since we will not rely on the results we
leave the corrections to the reader.} ~\cite{LQ}, p51, there exists
a further subsequence $(n_{i_j} =: n_j)$ and a control function $w$ such that
$$
     | X^k_{s,t} |, | X^k_{s,t} (n_j) | \le w(s,t)^{k/p}, \ \ \ k=1,2
$$
for any $j$ and
$$
   | X^k_{s,t} - X^k_{s,t} (n_j) | \le 2^{-j} w(s,t)^{k/p}, \ \ \ k=1,2.
$$
Hence, by the ULT estimates above,
$$
 | Y^k_{s,t} - Y^k_{s,t} (n_j) | \le 2^{-j} C w(s,t)^{k/p}, \ \ \ k=1,2.
$$
Using the super-additivity this implies $d(Y(n_j),Y) \to 0$. Since we were able
to extract, from any subsequence, a further convergent subsequence with limit $Y$,
it is clear that $d(Y(n),Y) \to 0$. Hence $X \to Y$ is continuous in $p$-variation metric. QED.
\\

{\bf Remarks: } (1) The fine-estimates exhibit some kind of Lipschitz behavior which is not visible in terms
of $p$-variation metric. \\
(2) Since a.s. $d(\beta(n),\beta) \to 0$ Theorem (\ref{ULT}) applies
and yields a Wong-Zakai-type result.
However, constructing a control $w$ as above is certainly a bad idea, since Corollary \ref{cor32} provides us with a much easier control, 
which indeed controls the whole sequence rather
than just a subsequence. Essentially, $w(s,t) = c|t-s|$ for some (random) constant  $c$.
We will now exploit these observations.

\section{Lipschitz regularity of the It\^o-map for H\"older Rough Paths}
 \label{sec5}

%
%
%

\begin{theorem} \label{thm52}

 Under assumptions on $f$ as in Theorem \ref{ULT} the It\^o-map $F: X \to Y$ is locally Lipschitz-continuous from $(H\Omega_p(V),\rho ) \to (H\Omega_p(W),\rho)$.
\end{theorem}

{\bf Proof: } Assume $ \rho (X , \hat X) \le \epsilon $. This just means that (always for all $s,t$)
$$
      | X^k_{s,t} - \hat X^k_{s,t} | \le \epsilon |t-s|^{k/p}, \ \ \ k =1,2
$$
Introduce an additive control function
$$
    w(s,t)  := ( \| X \| \lor \| \hat X \| \lor 1 )|t-s|.
$$
With this choice, the assumptions (\ref{ass1}) and (\ref{ass2}) are satisfied and Theorem \ref{ULT} 
tells us that there exists a constant $C$, depending on the maximum of the control $w$
and hence only on $\|X \| \lor \| \hat X \|$, such that

$$
    | Y^k_{s,t} - \hat Y^k_{s,t} | \le \epsilon (C w(s,t))^{k/p}, \ \ \ k=1,2.
$$
Expanding $w$ we obtain, for a new constant $C$ depending only on $\|X \| \lor \| \hat X \|$,
$$
   \rho(Y,\hat Y) = \| Y - \hat Y \| \le \epsilon C
$$ 

The conclusion follows.    QED

\section{Application to the Support Theorem}  \label{sec6}

Consider the Stratonovich-SDE

\begin{equation}
     dY = f(Y) \circ d\beta   \label{SDE}
\end{equation}

with initial condition $Y(0)=y_0 \in W=\Re^N$ where $f=(f_1,...,f_d)$ stands for $d$ vector-fields, $C^3$ and globally Lipschitz, driven by Brownian motion in $V = \Re^d$. (The case of an additional drift-term $f_0 (Y) dt$ is a trivial generalization of the method below.)
In Lyons' form this corresponds to the differential equation
$$
     dy = f (y) dx
$$
lifted to Rough Paths level and using Stratonovich enhanced BM $\beta$ as
driving signal. By Theorem \ref{thm52} the corresponding It\^o-map
$F: H\Omega_p (V) \to H\Omega_p (W)$ is continuous w.r.t. the metric $\rho$.
Since $F=(F^1,F^2)$ takes values in $W \oplus W^{\otimes 2}$ there is
a natural projection to its underlying path in $W$. Set
$$
          \Phi_t(\beta) = y_0 + [F(\beta)]^1_{0,t}.
$$
Note that $\Phi$ is continuous from $(H\Omega_p(V), \rho)$ to
$C^{0,1/p}([0,1],W)$ with usual H\"older semi-norm (actually norm,
since all paths are pinned at $y_0$ at time $0$).
Set
$$
           \Psi (\omega) = \Phi (\beta(\omega)).
$$
The map $\Psi$ maps $C([0,1],V)$ to $C^{0,1/p}([0,1],W)$ and is measurable
only (due to $\omega \to \beta(\omega)$) and we will also call it
It\^o-map (no confusion will arise). Note that
$$
         Y_t (\omega) := \Psi_t (\omega)
$$
solves (\ref{SDE}). Equip $C([0,1],V)$ with the standard
Wiener-measure $\mu$.  Our aim is to describe the support of 
$ (\Psi)_* \mu = (\Phi)_* P$ where $P$ is the law of EBM on $H\Omega_p(V)$.

Continuity of the $\Phi$ allows to restrict the
discussion to the support of $P$ only. 

\begin{theorem} The support of
$P$ equals the $\rho$-closure of $\mathcal D$. 
\end{theorem}

{\bf Proof: }


From Corollary \ref{cor34} resp. Corollary \ref{cor36}, a.s. and in $\rho$-metric,
$$
      \beta^{ad} (n) (\omega) \to \beta (\omega), \ \ \ 
      \beta(\omega - \omega^{ad}(n) +h) \to (h^1,h^2),
$$
for arbitrary $(h^1,h^2) \in {\mathcal D}$. The first convergence
implies that supp $P$ is contained in the $\rho$-closure of $\mathcal D$.
The second convergence result shows, in particular, convergence in probability
of $ \beta(\omega - \omega^{ad}(n) +h) \to (h^1,h^2)$. This means
that the probability of being within an $\epsilon-ball$ from $(h^1,h^2)$ tends 
to $1$, hence is positive for $n$ large enough. By Girsanov's theorem the
same is true for $\beta(\omega)$. Hence we get to other desired inclusdion,
${\mathcal D}$ is contained in supp $P$. (This argument is due to ~\cite{MS}).

{\bf Remark:} As in ~\cite{LQZ} the closure of ${\mathcal D}$ is seen to
coincide with the closure of the natural lift of the Cameron-Martin space. This also follows from the results in ~\cite{FV}.

Remarking that solving a rough differential equation driven by a smooth
rough path amounts to solve a controlled ODE (~\cite{LQ}, p164) we have

\begin{corollary}
The law of the diffusion-process $Y$ in equation (\ref{SDE}), that is,
$(\Psi)_* \mu$ is the ($1/p$-H\"older) closure of the solutions to the
 control ODE,
$$
                            dy=f(y)dh 
$$
for all $h \in {\mathcal D}$.
\end{corollary}

(Similarly, one can obtain a support description for the rough path solution
to (\ref{SDE}).)

\appendix

\section{Proof of Lemma (\ref{lemma21})}

{\bf Proof.} To avoid trivialities assume $s<t$, $s' < t'$. W.l.o.g. $s \le s'$. Then
\begin{eqnarray*}
 E[ | X_{s,t} - X_{s',t'} |^2 ] & = & E[ |X_{s,t}|^2] - 2 E [ X_{s,t} X_{s',t'} ] + E[|X_{s',t'}|^2] \\
      & = & \frac{1}{2} |t-s|^{2-2\alpha} - 2  E [ X_{s,t} X_{s',t'} ] + \frac{1}{2} |t'-s'|^{2-2\alpha}.
\end{eqnarray*}
To deal with the middle part we distinguish between a few cases. \\

Case i) $t \le s'$. By independence of Brownian increments, 

\begin{eqnarray*}
 E[ | X_{s,t} - X_{s',t'} |^2 ]
      & = & \frac{1}{2} \big( |t-s|^{2-2\alpha} + |t'-s'|^{2-2\alpha} \big) \\
      & \le & \frac{1}{2} \big( |s'-s|^{2-2\alpha} + |t'-t|^{2-2\alpha} \big) \\
      & \le & \Big( \frac{|s'-s|+|t'-t|}{2} \Big)^{2-2\alpha} \\
      & = & c (\alpha)  \Big( |s'-s| + |t' - t| \Big)^{2-2\alpha}
\end{eqnarray*}
where we used $\alpha \in [1/2,1)$ for the last inequality. \\


Case ii.1) $s \le s' < t \le t'$.
Since
$ 
    X_{s,t} - X_{s',t'} =  X_{s,t} - X_{s',t} + X_{s',t} - X_{s',t'}
$
we estimate seperatley
$
   E[ | X_{s,t} - X_{s',t}  |^2 ]
$
and
$
   E[ |  X_{s',t} - X_{s',t'} |^2 ].
$
As for the first, the cross-term is readily computed via Ito's isometry and we have
$$
 E[ | X_{s,t} - X_{s',t}  |^2 ] =  \frac{1}{2} |t-s|^{2-2\alpha} -
                                                \frac{ |t-s'|^{2-\alpha}}{|t-s|^\alpha }
						+  \frac{1}{2} |t-s'|^{2-2\alpha} .
$$
The r.h.s. is indeed bounded by a constant times $|s'-s|^{2-2\alpha}$ as follows
from the lemma below with $a =s'-s , b = t-s'$.
As for the second, we get
$$   E[ | X_{s',t} - X_{s',t'}  |^2 ] =  \frac{1}{2} |t-s'|^{2-2\alpha} -
                                                \frac{ |t'-s'|^{2-\alpha}}{|t-s'|^\alpha }
						+  \frac{1}{2} |t'-s'|^{2-2\alpha} .
$$
Choose $a=t'-t, b=t-s'$ and again apply the lemma. Both together yield the required estimate,
$$
  E[ |  X_{s,t} - X_{s',t'} |^2 ] \le c ( |s'-s| + |t'-t| )^{2-2\alpha}.
$$
Case ii.2) $s \le s' \le t' \le t$ Similar.   QED

\begin{lemma}
For non-negative real numbers $a,b$, $\alpha \in (0,1)$ one has
$$
    \frac{1}{2}(a+b)^{2-2\alpha} - \frac{ b^{2-\alpha} }{ (a+b)^\alpha } + \frac{1}{2}b^{2-2\alpha}  \le c a^{2-2\alpha}
$$
for some constant $c=c(\alpha)$.
\end{lemma}
{\bf Proof: } Divide the l.h.s. by $a^{2-2\alpha}$ and observe that it is a continuous function of $x=b/a \in [0,\infty)$.
An easy expansion shows that everything stays bounded as $x \to \infty$. The lemma follows. \\

{\bf Acknowledgement: } Thanks to S.Varadhan, G. Ben Arous and N. Victoir for related discussions, to E. Kosygina for reading an early draft. The author also
acknowledges gratefully financial support by the Austrian Academy of Science.

\end{document}